\documentclass[reqno,11pt]{amsart}

\usepackage[a4paper,textwidth=16.7cm,textheight=23.4cm,centering]{geometry}
\usepackage{amssymb,mathtools}
\usepackage{cite}
\usepackage[unicode,colorlinks=true,linkcolor=blue,citecolor=red,urlcolor=blue]{hyperref}
\usepackage{microtype}
\usepackage{needspace}

\numberwithin{equation}{section}
\mathtoolsset{showonlyrefs}
\allowdisplaybreaks[1]
\newtheorem{theorem}{Theorem}[section]
\newtheorem{corollary}[theorem]{Corollary}
\newtheorem{proposition}{Proposition}[section]
\newtheoremstyle{bracketnote}{3pt}{3pt}{\itshape}{}{\bfseries}{.}{.5em}{\thmname{#1}\thmnumber{ #2}\thmnote{ \normalfont#3}}
\theoremstyle{bracketnote}
\newtheorem{lemma}{Lemma}[section]
\theoremstyle{definition}

\theoremstyle{remark}
\newtheorem{remark}[theorem]{Remark}

\newcommand{\R}{\mathbb R}
\newcommand{\dd}{\mathop{}\!\mathrm{d}}

\title[Uniform velocity support for the repulsive Vlasov--Poisson system]
{\small Uniform Velocity Support for the Repulsive Vlasov--Poisson System}
\author{Zhaopeng Wang}
\thanks{School of Mathematics and Statistics, Beijing Institute of Technology, Beijing 100081, China.\newline Email: \texttt{zpwang@bit.edu.cn}}
\hypersetup{
  pdftitle={Uniform Velocity Support for the Repulsive Vlasov-Poisson System},
  pdfauthor={Zhaopeng Wang},
  pdfsubject={Uniform velocity support for the three-dimensional repulsive Vlasov-Poisson system},
  pdfkeywords={Vlasov-Poisson system; velocity support; dispersive estimates; Coulomb-corrected energy; convex transport; particle-passage estimates; Coulomb interaction; velocity moments}
}

\begin{document}

\begin{abstract}
This paper proves that the velocity support of the three-dimensional repulsive Vlasov--Poisson system remains uniformly bounded for all time for nonnegative initial data in $C_c^1(\mathbb R^6)$. The proof establishes uniform bounds on a weighted fifth-order energy moment and an associated fourth-order force integral by exploiting the relative motion of particle pairs. On intervals where the target particle moves sufficiently fast relative to the expansion, these estimates give a summable bound for the current in its energy identity. This bound yields uniform control of the velocity support.
\end{abstract}
\maketitle



\section{Introduction}\label{sec:introduction}

\subsection{The model and the velocity support}\label{subsec:model}

A plasma consists of a large number of charged particles. When collisions between particles can be neglected, their motion can be described by a kinetic equation. The unknown $f(t,x,v)$ is the distribution function in phase space. It describes the density of particles at position $x$ with velocity $v$ at time $t$.

The Vlasov--Poisson system is a classical model for the motion of a collisionless plasma under its self-consistent electric field. In three dimensions, the repulsive system is written as
\begin{equation}\label{eq:VP}
\left\{
\begin{aligned}
&\partial_t f+v\cdot\nabla_x f+E\cdot\nabla_v f=0,\quad f|_{t=0}=f_0,\quad t\geq0,\ x,v\in\mathbb R^3,\\
&E=-\nabla_x U,\quad -\Delta_x U=\rho,\quad \rho(t,x)=\int_{\mathbb R^3}f(t,x,v)\,\dd v.
\end{aligned}
\right.
\end{equation}
Here, $\rho(t,x)$ is the spatial density, $E(t,x)$ is the electric field, and $U(t,x)$ is the electric potential. With the condition $U(t,x)\to0$ as $|x|\to\infty$, the potential and the field are given by
\[
U(t,x)=\int_{\mathbb R^3}\frac{\rho(t,x')}{4\pi|x-x'|}\,\dd x',\qquad E(t,x)=\int_{\mathbb R^3}\frac{x-x'}{4\pi|x-x'|^3}\rho(t,x')\,\dd x'.
\]
With this sign, particles of the same charge repel each other.

The total energy of the system is defined by
\begin{equation}\label{eq:total-energy}
\mathcal E(t)=\frac12\int_{\mathbb R^6}|v|^2f(t,x,v)\,\dd x\,\dd v+\frac12\int_{\mathbb R^3}U(t,x)\rho(t,x)\,\dd x.
\end{equation}
The first term is the kinetic energy of the particles, and the second is the electrostatic potential energy generated by their self-consistent field. For classical solutions, the total energy in \eqref{eq:total-energy} is conserved: $\mathcal E(t)=\mathcal E(0)$.

For a particle starting from an initial point $(x_0,v_0)$, we denote its position and velocity at time $t$ by $X(t)$ and $V(t)$, respectively. They satisfy
\begin{equation}\label{eq:characteristics}
\dot X(t)=V(t),\qquad \dot V(t)=E(t,X(t)).
\end{equation}
Thus the particles generate the electric field, while the field acts back on their motion and changes their velocities.

This leads to a natural question. Even if all particles have bounded velocities initially, can the electric field accelerate some particles to arbitrarily large speeds as time goes on? To study this question, we define the radius of the velocity support by
\[
R_f(t)=\sup_{(x,v)\in\operatorname{supp}f(t)}|v|.
\]
The quantity $R_f(t)$ represents the largest speed among all particles present at time $t$. We are interested in whether the fastest particles in the cloud can be controlled uniformly over the whole evolution, that is, whether
\[
\sup_{t\geq0}R_f(t)<\infty.
\]

\subsection{Previous results and the main theorem}\label{subsec:previous-results}

The global classical theory of the three-dimensional Vlasov--Poisson system is now well established. Early global existence results include Batt's theorem for spherically symmetric systems~\cite{Batt1977} and the small-data result of Bardos and Degond~\cite{BardosDegond1985}. For general smooth and compactly supported initial data, the global existence and regularity theory was developed in the works of Lions and Perthame~\cite{LionsPerthame1991}, Pfaffelmoser~\cite{Pfaffelmoser1992}, and Schaeffer~\cite{Schaeffer1991}; see also the survey~\cite{Rein2007}. These results provide global classical solutions, with bounded velocity support on every finite time interval.

For the repulsive system, Illner and Rein~\cite{IllnerRein1996}, as well as Perthame~\cite{Perthame1996}, established time-decay estimates based on the dispersive effect of particle motion. These estimates describe the large-time expansion of the particle cloud and control the mean-square velocity relative to this expansion. They provide important information on the long-time behavior of general solutions.

For spherically symmetric solutions, Horst~\cite{Horst1990} established decay estimates for the spatial density and the electric field.  Pankavich~\cite{Pankavich2021} established uniform bounds on the velocity support for the classical Vlasov--Poisson system and on the momentum support for the relativistic system. He also obtained detailed information on the large-time behavior of the characteristics. For general, non-spherically symmetric data, whether the maximum particle speed admits a uniform bound for all time remains open. Rein also discussed the long-time boundedness of particle velocities in connection with a conjectured $t^{-3}$ decay of the spatial density in his survey~\cite{Rein2007}.

The large-time growth of the velocity support for general data has been studied by several authors. Rein obtained a sublinear growth bound with polynomial exponent $2/3$~\cite{Rein1998}. Pallard improved this estimate to $t^{2/7+\varepsilon}$ for any $\varepsilon>0$~\cite{Pallard2011a}, and later obtained a bound of order $(t\log t)^{6/25}$~\cite{Pallard2011b}. Further improvements were obtained by Chen and Zhang~\cite{ChenZhang2014}, who established the bound $t^{2/11+\varepsilon}$ for any $\varepsilon>0$, and then by Chen and Li~\cite{ChenLi2018}, who proved
\[
R_f(t)\leq C(1+t)^{2/15}\log^{8/15}(2+t).
\]
These estimates give sublinear upper bounds on the velocity support for general data.

In our recent work~\cite{WangGrowth}, we introduced a corrected energy associated with an individual particle trajectory. The main idea is to combine the kinetic energy relative to the linear expansion with the positive Coulomb potential evaluated along the same trajectory. In the resulting energy identity, the direct electric-field contributions cancel, leaving a scalar current generated by the source particles. This identity provides a way to connect estimates for the whole particle cloud with the motion of an individual particle, and it is the starting point of the present work.

In the present paper, we combine the corrected energy identity with estimates based on the relative motion of particle pairs. We first establish weighted energy and force bounds for the whole particle cloud. We then use these bounds to control the accumulated current along each individual particle trajectory. This yields a uniform bound on the velocity support for general compactly supported initial data, without any symmetry or smallness assumption.

Now we state our main result as follows.

\begin{theorem}\label{thm:main}
Let $0\leq f_0\in C_c^1(\mathbb R_x^3\times\mathbb R_v^3)$, and let $f$ be the corresponding global classical solution of \eqref{eq:VP}. There is a constant $C(f_0)<\infty$ such that
\[
\sup_{t\geq0}R_f(t)\leq C(f_0).
\]
\end{theorem}

\subsection{Expanding coordinates}\label{subsec:expanding-energy}

We compare the motion of each particle with a reference motion whose position grows linearly in time. This allows us to separate the large-scale expansion of the particle cloud from the motion of individual particles. Consider $x(t)=(1+t)y$, where $y$ is fixed. The corresponding velocity is $y=x/(1+t)$. Thus $x/(1+t)$ is the reference velocity at position $x$. We define
\begin{equation}\label{eq:expanding-coordinates}
y=\frac{x}{1+t},\qquad w=v-\frac{x}{1+t}.
\end{equation}
The variable $y$ records position on the expanding spatial scale, while $w$ is the velocity relative to the reference motion. The inverse transformation is
\begin{equation}\label{eq:physical-velocity}
x=(1+t)y,\qquad v=y+w.
\end{equation}
Time-dependent rescalings for the Vlasov--Poisson system were studied by Dolbeault and Rein~\cite[Section~3]{DolbeaultRein2001}. The spatial scaling above is also related to the projective transformations considered by Serre~\cite{SerreProjective}.

Along a characteristic satisfying \eqref{eq:characteristics}, differentiation of \eqref{eq:expanding-coordinates} gives
\[
\frac{\dd y}{\dd t}=\frac{V}{1+t}-\frac{X}{(1+t)^2}=\frac{w}{1+t},\qquad \frac{\dd w}{\dd t}=E(t,X)-\frac{w}{1+t}.
\]
To make the damping coefficient constant, introduce the logarithmic time $\tau=\log(1+t)$ and the rescaled field
\[
E^\sharp(\tau,y)=e^{2\tau}E(e^\tau-1,e^\tau y).
\]
Since $\dd t/\dd\tau=1+t$, the characteristic equations become
\begin{equation}\label{eq:expanding-characteristics}
\dot y=w,\qquad \dot w=-w+e^{-\tau}E^\sharp(\tau,y).
\end{equation}
From now on, a dot denotes differentiation with respect to $\tau$ along a particle trajectory. The term $-w$ is a linear damping term generated by the change of coordinates, not an additional force in the original Vlasov--Poisson system.

We also write the physical potential in the new coordinates and define the corrected energy by
\[
U^\sharp(\tau,y)=U(e^\tau-1,e^\tau y),\qquad H(\tau,y,w)=\frac12|w|^2+U^\sharp(\tau,y).
\]
The potential $U^\sharp$ is evaluated in the new coordinates without an additional scaling factor. Since $U\geq0$, the corrected energy controls the relative kinetic energy. Its evolution has a useful cancellation: the electric-field term in the derivative of $|w|^2/2$ cancels with the corresponding term in the derivative of the potential. Along a particle trajectory,
\[
\dot H+H+\frac12|w|^2=e^{-\tau}J^\sharp.
\]
The integral expression for the scalar current $J^\sharp$ is given in Section~\ref{subsec:basic-estimates}. This identity is the main tool for studying the motion of individual particles.

To prove Theorem~\ref{thm:main}, we establish the following decay estimate for the velocity relative to the linear expansion.

\begin{corollary}\label{cor:relative-decay}
Under the assumptions of Theorem~\ref{thm:main}, there is a constant $C(f_0)<\infty$ such that
\begin{equation}\label{eq:relative-decay}
\sup_{(x,v)\in\operatorname{supp}f(t)}\left|v-\frac{x}{1+t}\right|\leq C(f_0)\bigl(1+\log(1+t)\bigr)^{-3/2},\qquad t\geq0.
\end{equation}
\end{corollary}

\begin{remark}\label{rem:decay-rate}
The logarithmic decay rate in Corollary~\ref{cor:relative-decay} is sufficient for the present purpose. The power $3/2$ makes the relative velocity integrable in logarithmic time, which yields the uniform velocity-support bound in Theorem~\ref{thm:main}. We do not pursue a sharper decay rate here.
\end{remark}

Indeed, along every characteristic from $\operatorname{supp}f_0$, \eqref{eq:relative-decay} gives $|w(\tau)|\leq C(f_0)(1+\tau)^{-3/2}$. Using $\dot y=w$ from \eqref{eq:expanding-characteristics} and the compact initial spatial support, we obtain
\[
|y(\tau)|\leq |y(0)|+\int_0^\tau|w(\sigma)|\,\dd\sigma\leq |y(0)|+C(f_0)\int_0^\infty(1+\sigma)^{-3/2}\,\dd\sigma\leq C(f_0).
\]
The same decay estimate bounds $w$ uniformly. Since $V=y+w$ by \eqref{eq:physical-velocity} and the characteristic flow transports the initial support onto $\operatorname{supp}f(t)$, Corollary~\ref{cor:relative-decay} implies Theorem~\ref{thm:main}.

\subsection{Proof strategy}\label{subsec:proof-strategy}

A central difficulty in proving Theorem~\ref{thm:main} is to control the accumulated current along an individual particle trajectory. The corrected energy identity in Section~\ref{subsec:expanding-energy} reduces the evolution of a particle's energy to the current generated by the source particles. However, the available energy estimates concern the whole particle cloud, whereas the current in this identity acts along a single trajectory. Bridging these two levels of estimates is the main difficulty in the proof.

In the spherically symmetric setting, the electric field is radial and can be expressed in terms of the mass enclosed by a sphere centered at the origin. This geometric structure is useful in estimating the field along particle trajectories. For general data without spherical symmetry, the field receives contributions from particles in all directions, and this radial description is no longer available. We therefore need to control the particle interactions without relying on any symmetry of the solution.

Our approach is to use the relative motion of particle pairs to control the current over the whole particle cloud. We integrate the corrected energy identity with respect to the particle distribution and use the oddness of the Coulomb kernel to exchange the labels of two particles. In the kinetic contribution, this symmetrization produces a factor involving their relative velocity, which is also the velocity of their relative position. This allows us to apply the passage lemma to the singular Coulomb interaction. Together with the convex transport estimate, this provides the main tools for the proof. The convex transport estimate also allows us to use velocity weights. In particular, taking $\omega=|w|^3$ gives $\omega^{2/3}=|w|^2$, and hence a space-time $L^2$ estimate for the second velocity moment. This estimate is used in the quotient bound in Section~\ref{sec:energy-moments}.

The proof proceeds in two stages. In the first stage, we estimate the corrected energy over the whole particle cloud. The convex transport estimate controls the spatial density and the velocity-weighted density moments. The passage lemma then uses the relative motion of particle pairs to handle the singular Coulomb interaction. Combining these estimates, we establish weighted bounds on the $5/2$-moment of the corrected energy and on the fourth-order force integral. The rough density estimate controls the contribution from particles that are close to each other. These estimates are summarized in Proposition~\ref{prop:cloud-moment}.

In the second stage, we fix an arbitrary particle trajectory starting from the initial support and estimate the current acting on this particle. We use the estimates for the whole particle cloud obtained in Proposition~\ref{prop:cloud-moment}, together with the fixed-target version of the passage lemma. This gives a summable bound on the accumulated current over any time interval on which the relative speed of the target particle stays above a prescribed threshold. The estimate is uniform with respect to the choice of the target particle, and is stated in Proposition~\ref{prop:single-current}.

We then apply this estimate to the corrected energy identity along the fixed trajectory. A suitable time weight and a last-crossing argument show that the weighted particle energy remains uniformly bounded. The potential estimate ensures that, whenever this weighted energy exceeds a sufficiently large fixed level, the relative speed satisfies the threshold required in Proposition~\ref{prop:single-current}. Summing the current estimates over the corresponding time intervals yields the decay estimate in Corollary~\ref{cor:relative-decay}. Since this decay is integrable in logarithmic time, the particle position in the expanding coordinates remains bounded. Together with the bound on the relative velocity, this gives the uniform velocity-support estimate in Theorem~\ref{thm:main}.

The rest of the paper is organized as follows. Section \ref{sec:preparations} collects the basic estimates and proves the convex transport estimate and the passage lemma. Section \ref{sec:energy-moments} establishes the weighted energy-moment and force bounds for the whole particle cloud. Section \ref{sec:proof-main} uses these bounds to control the current along a fixed particle trajectory and completes the proof of Theorem~\ref{thm:main}.

\section{Preparations}\label{sec:preparations}

\subsection{Basic estimates}\label{subsec:basic-estimates}

We use the expanding coordinates introduced in Section~\ref{subsec:expanding-energy}. Recall that $x=(1+t)y$, $v=y+w$, and $\tau=\log(1+t)$. The Jacobian of this change of variables is
\[
\left|\frac{\partial(x,v)}{\partial(y,w)}\right|=(1+t)^3.
\]
We therefore define the distribution function in the new coordinates by
\[
g(\tau,y,w)=(1+t)^3f(t,(1+t)y,y+w),\qquad t=e^\tau-1.
\]
In particular, $g(\tau,y,w)\,dy\,dw=f(t,x,v)\,dx\,dv$. Thus the transformation preserves the total particle mass. Since $f$ is constant along the original characteristics, we also have $\|g(\tau)\|_{L^\infty_{y,w}}=e^{3\tau}\|f_0\|_{L^\infty_{x,v}}$.

For $k\geq0$, define the spatial velocity moments and their total integrals by
\[
m_k(\tau,y)=\int_{\mathbb R^3}|w|^kg(\tau,y,w)\,dw,\qquad M_k(\tau)=\int_{\mathbb R^3}m_k(\tau,y)\,dy.
\]
In particular, $M_0$ is the conserved total mass, and $m_0$ is the spatial density in the expanding coordinates. Write $K(z)=z/(4\pi|z|^3)$. By the homogeneity of this kernel, the rescaled electric field satisfies
\[
E^\sharp(\tau,y)=\int_{\mathbb R^3}K(y-y')m_0(\tau,y')\,dy'.
\]

The characteristic equations in the expanding coordinates are $\dot y=w$ and $\dot w=-w+e^{-\tau}E^\sharp(\tau,y)$. Consequently, $g$ satisfies the phase-space continuity equation
\begin{equation}\label{eq:transformed-continuity}
\partial_\tau g+\nabla_y\cdot(wg)+\nabla_w\cdot\bigl((e^{-\tau}E^\sharp-w)g\bigr)=0.
\end{equation}
For a sufficiently regular function $\psi(\tau,y,w)$, integration by parts gives
\begin{equation}\label{eq:phase-space-integration}
\frac{d}{d\tau}\iint_{\mathbb R^6}\psi g\,dy\,dw=\iint_{\mathbb R^6}\left[\partial_\tau\psi+w\cdot\nabla_y\psi+(e^{-\tau}E^\sharp-w)\cdot\nabla_w\psi\right]g\,dy\,dw.
\end{equation}
The expression in brackets is the derivative of $\psi$ along a characteristic; we refer to \eqref{eq:phase-space-integration} below. Taking $\psi=1$ shows again that the total mass is conserved. Taking functions depending only on $y$ gives the spatial continuity equation
\begin{equation}\label{eq:spatial-continuity}
\partial_\tau m_0+\nabla_y\cdot j=0,\qquad j(\tau,y)=\int_{\mathbb R^3}wg(\tau,y,w)\,dw.
\end{equation}

We next recall the energy identity along one characteristic and give its current term explicitly. For a particle with label $i$, let $H_i(\tau)=\frac12|w_i(\tau)|^2+U^\sharp(\tau,y_i(\tau))$. The kinetic energy satisfies
\[
\frac{d}{d\tau}\frac12|w_i|^2+|w_i|^2=e^{-\tau}w_i\cdot E^\sharp(\tau,y_i).
\]
On the other hand, differentiating the Coulomb potential along the same trajectory gives
\[
\frac{d}{d\tau}U^\sharp(\tau,y_i)+U^\sharp(\tau,y_i)=e^{-\tau}J^\sharp(\tau,y_i)-e^{-\tau}w_i\cdot E^\sharp(\tau,y_i),
\]
where
\[
J^\sharp(\tau,y_i)=\iint_{\mathbb R^6}K(y_i-y_j)\cdot w_j\,g(\tau,y_j,w_j)\,dy_j\,dw_j.
\]
The electric-field terms cancel because the kinetic and potential energies are differentiated along the same particle trajectory. We obtain
\begin{equation}\label{eq:single-energy-2}
\frac{d}{d\tau}H_i+H_i+\frac12|w_i|^2=e^{-\tau}J^\sharp(\tau,y_i).
\end{equation}
 
We now turn to the energy of the whole cloud. Define the kinetic and potential energies by
\[
\mathcal K(\tau)=\frac12M_2(\tau),\qquad \mathcal P(\tau)=\frac12\int_{\mathbb R^3}U^\sharp(\tau,y)m_0(\tau,y)\,dy.
\]
The potential energy also has the representation
\[
\mathcal P(\tau)
=\frac{e^{-\tau}}{2}
\iint_{\mathbb R^3\times\mathbb R^3}
\frac{m_0(\tau,y_i)m_0(\tau,y_j)}{4\pi|y_i-y_j|}
\,dy_i\,dy_j.
\]
The factor $1/2$ accounts for the two orders in which each pair is counted,
while $e^{-\tau}$ comes from the homogeneity of the Coulomb potential. 

\begin{lemma}\label{lem:basic-energy}
There exists a constant $C(f_0)>0$ such that, for every $\tau\geq0$,
\begin{equation}\label{eq:basic-moment-energy}
\mathcal K(\tau)+\mathcal P(\tau)\leq C(f_0)e^{-\tau}.
\end{equation}
Moreover,
\begin{align}
\|E^\sharp(\tau)\|_{L^2_y}&\leq C(f_0),\label{eq:field-l2}\\
\|U^\sharp(\tau)\|_{L^\infty_y}&\leq C(f_0)e^{-\tau/2}.\label{eq:potential-decay-basic}
\end{align}
\end{lemma}

\begin{proof}
Multiplying \eqref{eq:transformed-continuity} by $|w|^2/2$ and integrating over phase space gives
\[
\frac{d}{d\tau}\mathcal K+2\mathcal K=e^{-\tau}\int_{\mathbb R^3}j\cdot E^\sharp\,dy.
\]
On the other hand, differentiating the double-integral representation of $\mathcal P$ and using \eqref{eq:spatial-continuity} yields
\[
\frac{d}{d\tau}\mathcal P+\mathcal P=-e^{-\tau}\int_{\mathbb R^3}j\cdot E^\sharp\,dy.
\]
The field terms cancel. Hence
\[
\frac{d}{d\tau}(\mathcal K+\mathcal P)+2\mathcal K+\mathcal P=0.
\]
Since both energies are nonnegative, integration gives \eqref{eq:basic-moment-energy}.

The Poisson equation and integration by parts imply
\[
\|E^\sharp(\tau)\|_{L^2_y}^2=(1+t)\|E(t)\|_{L^2_x}^2=(1+t)\int_{\mathbb R^3}U(t,x)\rho(t,x)\,dx=2e^\tau\mathcal P(\tau)\leq C(f_0).
\]
This proves the estimate for the electric field in \eqref{eq:field-l2}.

We next derive a pointwise bound on the potential. For fixed $x$, split the velocity integral defining $\rho(t,x)$ into the regions $|v-x/(1+t)|\leq h$ and $|v-x/(1+t)|>h$. The first part is bounded by $C\|f_0\|_\infty h^3$, while the second is bounded by
\[
h^{-2}\int_{\mathbb R^3}\left|v-\frac{x}{1+t}\right|^2f(t,x,v)\,dv.
\]
Optimizing in $h$ gives
\[
\rho(t,x)\leq C\|f_0\|_\infty^{2/5}\left(\int_{\mathbb R^3}\left|v-\frac{x}{1+t}\right|^2f(t,x,v)\,dv\right)^{3/5}.
\]
Raising this estimate to the power $5/3$, integrating in $x$, and using \eqref{eq:basic-moment-energy}, we obtain
\[
\|\rho(t)\|_{L^{5/3}_x}\leq C(f_0)(1+t)^{-3/5}.
\]
Splitting the Newton potential into a ball of radius $r$ and its complement, Hölder's inequality and conservation of mass give
\[
\|U(t)\|_{L^\infty_x}\leq C\left(r^{1/5}\|\rho(t)\|_{L^{5/3}_x}+r^{-1}M_0\right).
\]
Taking $r=(1+t)^{1/2}$ and recalling the definition of $U^\sharp$, we obtain
\[
\|U^\sharp(\tau)\|_{L^\infty_y}=\|U(t)\|_{L^\infty_x}\leq C(f_0)e^{-\tau/2}.
\]
This proves \eqref{eq:potential-decay-basic} and completes the proof. \end{proof}

We also need a rough bound on the support of $g$ and its spatial density. We record these estimates in the following lemma.

\begin{lemma}\label{lem:rough-support}
There exists a constant $C(f_0)>0$ such that, for every $t\geq0$, with $\tau=\log(1+t)$,
\begin{equation}\label{eq:rough-velocity}
R_f(t)\leq C(f_0)(1+t),
\end{equation}
and, in the expanding coordinates,
\begin{equation}\label{eq:rough-support}
|y|+|w|\leq C(f_0)e^\tau\quad\text{on }\operatorname{supp}g(\tau),
\end{equation}
\begin{equation}\label{eq:rough-density}
\|m_0(\tau)\|_{L^\infty_y}\leq C(f_0)e^{6\tau}.
\end{equation}
\end{lemma}

\begin{proof}
The first estimate is the weaker velocity-support bound that follows from Chen and Li~\cite{ChenLi2018}. Since the initial spatial support is compact and $\dot X=V$ along characteristics, it gives $|X(t)|\leq C(f_0)(1+t)^2$. By the definition of the expanding coordinates, $y=X/(1+t)$ and $w=V-X/(1+t)$. Hence $|y|+|w|\leq C(f_0)(1+t)=C(f_0)e^\tau$, proving \eqref{eq:rough-support}.

On the other hand, the $L^\infty$ estimate for $g$ and the velocity support bound imply
\[
\begin{aligned}
m_0(\tau,y)&=\int_{\mathbb R^3}g(\tau,y,w)\,dw\\
&\leq \|g(\tau)\|_{L^\infty_{y,w}}\,|\operatorname{supp}_w g(\tau,y,\cdot)|\leq C(f_0)e^{3\tau}e^{3\tau}.
\end{aligned}
\]
Taking the supremum in $y$ proves \eqref{eq:rough-density} and completes the proof. \end{proof}

\subsection{Convex transport and particle passages}\label{subsec:two-lemmas}

The estimates in this subsection provide two complementary tools for the arguments in the next sections. The first controls spatial densities through convex transport, while the second controls the singular Coulomb interaction by using the relative motion of particle pairs. Together, they allow us to pass from estimates on the whole particle cloud to estimates along an individual particle trajectory.

The first lemma is a convex transport estimate for spatial densities, based on Brenier's convex transport construction~\cite{Brenier1991} and Serre's compensated integrability method~\cite{Serre2019}. Related ideas were also used in our previous work on moment propagation and energy conservation for the plasma--charge model~\cite{WangPlasmaCharge}. Here we adapt this method to the expanding coordinates and to velocity-weighted spatial densities.

For an integer $n\geq0$, an interval $J\subset[n,n+1]$, and a positive number $\ell$, define
\[
V_\ell(\omega)=\ell\sup_\Phi\int_J\iint_{\R^6}\omega g\,D^2\Phi(\ell\tau,y)[(1,w/\ell),(1,w/\ell)]\,\dd y\,\dd w\,\dd\tau,
\]
where the supremum is taken over all smooth convex functions $\Phi$ on $\R\times\R^3$ satisfying $|\nabla\Phi|\leq1$. The Hessian is taken in all four variables.

\begin{lemma}\label{lem:convex-transport}
For every nonnegative measurable weight $\omega$ that is bounded on the particle support over $J$,
\begin{equation}\label{eq:convex-density}
\left\|\int_{\mathbb R^3}\omega^{2/3}g\,dw\right\|_{L^2(J\times\mathbb R^3)}\lesssim e^n\ell^{1/2}V_\ell(\omega)^{2/3}.
\end{equation}
If $\omega$ is absolutely continuous along almost every characteristic, then
\begin{equation}\label{eq:convex-variation}
V_\ell(\omega)\leq 2\sup_{\tau\in J}\iint_{\R^6}\omega(1+|w|/\ell)g\,\dd y\,\dd w+\int_J\iint_{\R^6}\left(|\dot\omega|+\ell^{-1}\left|\frac{\dd}{\dd\tau}(w\omega)\right|\right)g\,\dd y\,\dd w\,\dd\tau.
\end{equation}
 
\end{lemma}

\begin{proof}
Fix an interval $J\subset[n,n+1]$ and a number $\ell>0$. We first change the time and velocity variables by setting $\sigma=\ell\tau$ and $\nu=w/\ell$. Write $X=(\sigma,y)\in\R^4$, and define
\[
\widetilde g(X,\nu)=\ell^3g(\sigma/\ell,y,\ell\nu),
\qquad
\widetilde\omega(X,\nu)=\omega(\sigma/\ell,y,\ell\nu).
\]
Set $\widetilde g=0$ outside $\sigma\in\ell J$, extend $\widetilde\omega$ by zero off the particle support, and define
\[
R(X)=\int_{\R^3}\widetilde\omega^{2/3}\widetilde g\,\dd\nu,
\qquad
D=\int_{\R^4}R(X)^2\,\dd X.
\]
On the fixed interval $J$, the classical solution has bounded phase-space support. The boundedness of $\omega$ therefore implies that $R$ is bounded and compactly supported, so $D<\infty$. We may assume that $D>0$ and $V_\ell(\omega)<\infty$, since the other cases are immediate.

The same change of variables rewrites the quantity in the lemma as
\[
V_\ell(\omega)=\sup_\Phi\int_{\R^4}\int_{\R^3}
\widetilde\omega\widetilde g\,D^2\Phi(X)[(1,\nu),(1,\nu)]\,\dd\nu\,\dd X,
\]
where the supremum is taken over smooth convex functions on $\R^4$ satisfying $|\nabla\Phi|\le1$.

Let $B_1^4$ be the unit ball in $\R^4$. Brenier's theorem~\cite{Brenier1991} provides a convex potential $\Phi_0$ whose gradient transports the probability density $R^2/D$ to the uniform probability density on $B_1^4$:
\[
(\nabla\Phi_0)_\#
\left(\frac{R(X)^2}{D}\,\dd X\right)
=
\frac{\mathbf 1_{B_1^4}(\zeta)}{|B_1^4|}\,\dd\zeta.
\]

We now construct a globally $1$-Lipschitz convex potential that induces the same transport. Let $S$ be a set of full source measure on which $\Phi_0$ is differentiable and $|\nabla\Phi_0|\le1$. Define
\[
\Phi(X)=\sup_{Z\in S}
\bigl\{\Phi_0(Z)+\nabla\Phi_0(Z)\cdot(X-Z)\bigr\}.
\]
Fix $X_0\in S$. The supporting-plane inequality and the bound on the slopes give
\[
\Phi_0(X_0)-|X-X_0|\le\Phi(X)\le\Phi_0(X_0)+|X-X_0|.
\]
Thus $\Phi$ is finite everywhere, convex, and globally $1$-Lipschitz.

For $X\in S$, we have $\Phi(X)=\Phi_0(X)$, and $\nabla\Phi_0(X)$ is a subgradient of $\Phi$ at $X$. Since $\Phi$ is differentiable almost everywhere and the source probability is absolutely continuous, $\nabla\Phi=\nabla\Phi_0$ almost everywhere with respect to that probability. Hence $\Phi$ induces the same transport.

Let $A(X)$ denote the Alexandrov Hessian of $\Phi$. The almost-everywhere Jacobian equation~\cite[Theorem~4.4 and Remark~4.5]{McCann1997} gives
\[
\det A(X)=\frac{|B_1^4|}{D}R(X)^2
\qquad\text{for almost every }X\in\{R>0\}.
\]
In particular, $A(X)$ is positive definite almost everywhere on $\{R>0\}$.

For any positive definite symmetric $4\times4$ matrix $M$, completing the square in $\nu$ and making a linear change of variables give
\[
\int_{\R^3}\bigl(M[(1,\nu),(1,\nu)]\bigr)^{-2}\,\dd\nu
=\frac{\pi^2}{\sqrt{\det M}}.
\]
At each fixed $X$, Hölder's inequality therefore yields
\[
\begin{aligned}
R(X)
&\le
\left(\int_{\R^3}\widetilde\omega\widetilde g\,
M[(1,\nu),(1,\nu)]\,\dd\nu\right)^{2/3}
\left(\int_{\R^3}
\frac{\widetilde g}{\bigl(M[(1,\nu),(1,\nu)]\bigr)^2}\,\dd\nu\right)^{1/3}\\
&\le
\pi^{2/3}\|\widetilde g\|_\infty^{1/3}(\det M)^{-1/6}
\left(\int_{\R^3}\widetilde\omega\widetilde g\,
M[(1,\nu),(1,\nu)]\,\dd\nu\right)^{2/3}.
\end{aligned}
\]
Raising this inequality to the power $3/2$ and taking $M=A(X)$ on $\{R>0\}$, we obtain
\[
|B_1^4|^{1/4}D^{-1/4}R(X)^2
\le
\pi\|\widetilde g\|_\infty^{1/2}
\int_{\R^3}\widetilde\omega\widetilde g\,
A(X)[(1,\nu),(1,\nu)]\,\dd\nu.
\]
The inequality is also valid almost everywhere on $\{R=0\}$, since its left-hand side vanishes and its right-hand side is nonnegative.

To bound the last integral by $V_\ell(\omega)$, choose a nonnegative $\eta\in C_c^\infty(B_1^4)$ with integral one, and set $\eta_\varepsilon(X)=\varepsilon^{-4}\eta(X/\varepsilon)$ and $\Phi_\varepsilon=\Phi*\eta_\varepsilon$. Each $\Phi_\varepsilon$ is smooth and convex, and
\[
|\nabla\Phi_\varepsilon(X)|
\le\int_{\R^4}|\nabla\Phi(X-h)|\eta_\varepsilon(h)\,\dd h\le1.
\]
Thus $\Phi_\varepsilon$ is admissible in the definition of $V_\ell(\omega)$.

At every Alexandrov point, the second-order expansion reads
\[
\Phi(X+h)=\Phi(X)+\nabla\Phi(X)\cdot h
+\tfrac12h^{\mathsf T}A(X)h+r_X(h),
\qquad r_X(h)=o(|h|^2).
\]
Insert this expansion into
\[
D^2\Phi_\varepsilon(X)
=\int_{\R^4}\Phi(X-h)D^2\eta_\varepsilon(h)\,\dd h.
\]
The constant and linear terms vanish, and the quadratic term gives $A(X)$. The remainder satisfies
\[
\bigl|D^2\Phi_\varepsilon(X)-A(X)\bigr|
\le
\sup_{0<|h|\le\varepsilon}\frac{|r_X(-h)|}{|h|^2}
\int_{\R^4}|z|^2|D^2\eta(z)|\,\dd z
\longrightarrow0.
\]
All the weighted quadratic forms are nonnegative. Fatou's lemma and the definition of $V_\ell(\omega)$ give
\[
\begin{aligned}
&\int_{\R^4}\int_{\R^3}\widetilde\omega\widetilde g\,
A(X)[(1,\nu),(1,\nu)]\,\dd\nu\,\dd X\\
&\quad\le
\liminf_{\varepsilon\downarrow0}
\int_{\R^4}\int_{\R^3}\widetilde\omega\widetilde g\,
D^2\Phi_\varepsilon(X)[(1,\nu),(1,\nu)]\,\dd\nu\,\dd X
\le V_\ell(\omega).
\end{aligned}
\]
Integrating the pointwise density estimate, we find
\[
|B_1^4|^{1/4}D^{3/4}
\le\pi\|\widetilde g\|_\infty^{1/2}V_\ell(\omega),
\qquad
D^{1/2}\le C\|\widetilde g\|_\infty^{1/3}V_\ell(\omega)^{2/3}.
\]
Finally, the change of variables and the bound $\|g(\tau)\|_\infty=e^{3\tau}\|f_0\|_\infty$ imply
\[
D^{1/2}
=\ell^{1/2}\left\|\int_{\R^3}\omega^{2/3}g\,\dd w\right\|_{L^2(J\times\R^3)},
\qquad
\|\widetilde g\|_\infty^{1/3}
\le\ell e^{n+1}\|f_0\|_\infty^{1/3}.
\]
Substitution and division by $\ell^{1/2}$ yield
\[
\left\|\int_{\R^3}\omega^{2/3}g\,\dd w\right\|_{L^2(J\times\R^3)}
\le C e^n\ell^{1/2}V_\ell(\omega)^{2/3}.
\]

Assume now that $\omega$ is absolutely continuous along almost every characteristic. If the right-hand side of \eqref{eq:convex-variation} is infinite, there is nothing to establish, so we assume it is finite. Fix a smooth convex test function $\Phi$ with $|\nabla\Phi|\le1$. Along each characteristic, $\dot y=w$, and hence
\[
\frac{\dd}{\dd\tau}(\ell\tau,y(\tau))=\ell(1,w/\ell).
\]
The product rule therefore gives
\[
\begin{aligned}
\frac{\dd}{\dd\tau}
\left[\omega\left(\partial_\sigma\Phi+\frac{w}{\ell}\cdot\nabla_y\Phi\right)\right]
=\ell\omega\,D^2\Phi[(1,w/\ell),(1,w/\ell)]+\dot\omega\,\partial_\sigma\Phi
+\frac1\ell\frac{\dd}{\dd\tau}(w\omega)\cdot\nabla_y\Phi,
\end{aligned}
\]
where the derivatives of $\Phi$ are evaluated at $(\ell\tau,y)$.

Integrate the preceding identity along almost every characteristic from $\tau_0$ to $\tau_1$, using the endpoint values of the absolutely continuous representative of $\omega$. Then integrate in the initial particle labels; the assumed integrability justifies this step. At time $\tau$, the characteristic flow carries the initial particle distribution to $g(\tau,y,w)\,\dd y\,\dd w$. On every compact subinterval $[\tau_0,\tau_1]\subset J$, we obtain
\[
\begin{aligned}
&\ell\int_{\tau_0}^{\tau_1}\iint_{\R^6}
\omega g\,D^2\Phi[(1,w/\ell),(1,w/\ell)]\,\dd y\,\dd w\,\dd\tau\\
&\quad=
\left[\iint_{\R^6}
\omega\left(\partial_\sigma\Phi+\frac{w}{\ell}\cdot\nabla_y\Phi\right)
g\,\dd y\,\dd w\right]_{\tau_0}^{\tau_1}-
\int_{\tau_0}^{\tau_1}\iint_{\R^6}
\left(\dot\omega\,\partial_\sigma\Phi
+\frac1\ell\frac{\dd}{\dd\tau}(w\omega)\cdot\nabla_y\Phi\right)
g\,\dd y\,\dd w\,\dd\tau.
\end{aligned}
\]

Since $|\partial_\sigma\Phi|\le1$ and $|\nabla_y\Phi|\le1$, the two endpoint contributions are bounded by $2\sup_{\tau\in J}\iint\omega(1+|w|/\ell)g\,\dd y\,\dd w$. Bounding the remaining terms by their absolute values, exhausting $J$ by compact subintervals, and taking the supremum over $\Phi$, we obtain
\[
\begin{aligned}
V_\ell(\omega)
&\le2\sup_{\tau\in J}\iint_{\R^6}
\omega(1+|w|/\ell)g\,\dd y\,\dd w+\int_J\iint_{\R^6}
\left(|\dot\omega|+\frac1\ell\left|\frac{\dd}{\dd\tau}(w\omega)\right|\right)
g\,\dd y\,\dd w\,\dd\tau.
\end{aligned}\qedhere
\]
\end{proof}

The second lemma uses the relative motion of two particles to control the Coulomb interaction accumulated during a close passage. Let $(y_i,w_i)$ and $(y_j,w_j)$ be two particle trajectories, and set $z=y_i-y_j$ and $\xi=w_i-w_j$. By the characteristic equations, $\dot z=\xi$.

Let $I_n=[n,n+1]$, where $n=0,1,2,\ldots$, and let $J\subset I_n$ be an interval. By Lemma~\ref{lem:rough-support}, we may choose $C_0\geq1$, depending only on the initial data, such that $|z|\leq C_0e^n$ on $I_n$. Set
\[
h=e^{-14n},\qquad k_h(z)=\frac{\mathbf 1_{\{|z|\leq C_0e^n\}}}{(h+|z|)^3}.
\]
Write $g_i=g(\tau,y_i,w_i)$ and $g_j=g(\tau,y_j,w_j)$. For a nonnegative measurable weight $B$ of a particle pair, write $dY=dy_i\,dw_i\,dy_j\,dw_j$ and define
{\small
\[
C_B:=2\sup_{\tau\in J}\int_{\mathbb R^{12}}|B\xi|g_i g_j\,dY+\int_J\int_{\mathbb R^{12}}\left|\frac{d}{d\tau}(B\xi)\right|g_i g_j\,dY\,d\tau.
\]
}

\begin{lemma}\label{lem:passage}
Suppose that $B\xi$ is absolutely continuous along almost every pair of particle trajectories and that $C_B<\infty$. Let $W$ be a measurable weight satisfying $0\leq W\leq C_1e^{5n}$ and $W=0$ almost everywhere on $\{B|\xi|^2=0\}$ with respect to the measure $g_i g_j\,dy_i\,dw_i\,dy_j\,dw_j\,d\tau$. Define $W^2/(B|\xi|^2)$ to be zero on this set. Then
\begin{equation}\label{eq:passage}
\begin{aligned}
&\int_J\iint_{\mathbb R^{12}}W|K(z)|g_i g_j\,dy_i\,dw_i\,dy_j\,dw_j\,d\tau\\
&\quad\lesssim e^{-3n}+\sqrt{(1+n)C_B}
\left(\int_J\iint_{\mathbb R^{12}}k_h(z)\frac{W^2}{B|\xi|^2}g_i g_j\,dy_i\,dw_i\,dy_j\,dw_j\,d\tau\right)^{1/2}.
\end{aligned}
\end{equation}
\end{lemma}

\begin{proof}
We first prove the estimate with integration over both particle labels. For $s\in\mathbb R$, define
\[
\psi_h(s)=\int_0^s\frac{dr}{h+|r|}.
\]
Then $\psi_h'(s)=(h+|s|)^{-1}$ and $|\psi_h(s)|=\log(1+|s|/h)$. Since $|z|\leq C_0e^n$ on the particle support and $h=e^{-14n}$, we have
\[
|\psi_h(z_k)|\leq\log(1+C_0e^{15n})\lesssim1+n,\qquad k=1,2,3.
\]
Using $\dot z_k=\xi_k$, the product rule gives, along almost every pair of characteristics,
\[
\frac{d}{d\tau}\big[(B\xi_k)\psi_h(z_k)\big]=(B\xi_k)'\psi_h(z_k)+\frac{B\xi_k^2}{h+|z_k|}.
\]
Fix a compact subinterval $[\tau_0,\tau_1]\subset J$. Integrating this identity along a pair of characteristics, with endpoint values taken from the absolutely continuous representative of $B\xi$, yields
\[
\int_{\tau_0}^{\tau_1}\frac{B\xi_k^2}{h+|z_k|}\,d\tau\lesssim(1+n)\left(|B\xi(\tau_0)|+|B\xi(\tau_1)|+\int_{\tau_0}^{\tau_1}|(B\xi)'|\,d\tau\right).
\]
We then integrate in the initial particle labels. The characteristic flow carries the initial particle distribution to $g(\tau,y,w)\,dy\,dw$. Hence the endpoint terms are controlled by the first term in $C_B$, while the time integral is controlled by its second term. Since $|z_k|\leq|z|$, we have
\[
\frac{B|\xi|^2}{h+|z|}\leq\sum_{k=1}^3\frac{B\xi_k^2}{h+|z_k|}.
\]
Letting the compact subintervals increase to $J$,  we obtain
\begin{equation*}
\int_J\iint_{\mathbb R^{12}}\frac{B|\xi|^2}{h+|z|}g_i g_j\,dy_i\,dw_i\,dy_j\,dw_j\,d\tau\lesssim(1+n)C_B.
\end{equation*}

We next consider the region $|z|\geq h$. Since $|K(z)|=(4\pi|z|^2)^{-1}\lesssim(h+|z|)^{-2}$, we use the factorization
\[
\frac{W}{(h+|z|)^2}=\left(\frac{B|\xi|^2}{h+|z|}\right)^{1/2}\left(\frac{W^2}{B|\xi|^2(h+|z|)^3}\right)^{1/2}.
\]
The identity holds almost everywhere by the convention in the statement. Applying Cauchy--Schwarz in time and in the particle labels, and using the preceding estimate, gives
\[
\begin{aligned}
&\int_J\int_{\{|z|\geq h\}}W|K(z)|g_i g_j\,dy_i\,dw_i\,dy_j\,dw_j\,d\tau\\
&\quad\lesssim\sqrt{(1+n)C_B}\left(\int_J\iint_{\mathbb R^{12}}k_h(z)\frac{W^2}{B|\xi|^2}g_i g_j\,dy_i\,dw_i\,dy_j\,dw_j\,d\tau\right)^{1/2}.
\end{aligned}
\]
The cutoff in $k_h$ removes no contribution, since $|z|\leq C_0e^n$ on the particle support.

It remains to estimate the region $|z|<h$. By Lemma~\ref{lem:rough-support},
\[
\|m_0(\tau)\|_{L^\infty_y}\lesssim e^{6n},\qquad \tau\in J.
\]
Using $W\leq C_1e^{5n}$, for each fixed target particle we find
\[
\begin{aligned}
&\iint_{|y_i-y_j|<h}W|K(y_i-y_j)|g_j\,dy_j\,dw_j\\
&\quad\leq C_1e^{5n}\|m_0(\tau)\|_{L^\infty_y}\int_{|z|<h}\frac{dz}{4\pi|z|^2}\lesssim e^{11n}h=e^{-3n}.
\end{aligned}
\]
Here $\int_{|z|<h}(4\pi|z|^2)^{-1}\,dz=h$. Integration in time contributes at most $|J|\leq1$, and integration in the target particle label contributes only the conserved total mass. Thus the near-field contribution is bounded by $Ce^{-3n}$. Combining the near-field and far-field estimates proves \eqref{eq:passage} for the double-particle integral.

\end{proof}

\begin{remark}\label{rem:passage-fixed-target}
The same estimate holds with one target particle $i$ fixed along any prescribed characteristic issuing from an initial label in $\operatorname{supp}f_0$, with integration only over the source particle $j$. The estimate is uniform in the choice of the target label; hence it applies to every characteristic from the initial support, including those starting at its boundary. In this case, the phase-space integrals in \eqref{eq:passage} are taken only over $(y_j,w_j)$ with measure $g_j\,dy_j\,dw_j$. The assumptions on $B\xi$ and $W$ are imposed along almost every source trajectory, and the zero-set condition is understood with respect to $g_j\,dy_j\,dw_j\,d\tau$. We write
{\small
\begin{equation}\label{eq:passage-fixed-variation}
C_{B,i}:=2\sup_{\tau\in J}\int_{\mathbb R^6}|B\xi|g_j\,dy_j\,dw_j+\int_J\int_{\mathbb R^6}\left|\frac{d}{d\tau}(B\xi)\right|g_j\,dy_j\,dw_j\,d\tau<\infty.
\end{equation}
}
In the fixed-target estimate, $C_B$ in \eqref{eq:passage} is replaced by $C_{B,i}$. All derivatives are taken along both the target and source trajectories. The implicit constant depends only on $f_0$ and $C_1$, and is independent of $n$, $J$, and the selected target trajectory.

Indeed, fix a target characteristic $i$. For almost every source characteristic $j$, the same product identity used in the proof of Lemma~\ref{lem:passage} holds, with $z=y_i-y_j$ and $\xi=w_i-w_j$. Integrating it from $\tau_0$ to $\tau_1$ and then only in the initial source label gives
\[
\begin{aligned}
\int_{\tau_0}^{\tau_1}\int_{\mathbb R^6}\frac{B|\xi|^2}{h+|z|}g_j\,dy_j\,dw_j\,&d\tau
\lesssim{}(1+n)\Bigg(\int_{\mathbb R^6}|B(\tau_0)\xi(\tau_0)|g(\tau_0,y_j,w_j)\,dy_j\,dw_j \\&+\int_{\mathbb R^6}|B(\tau_1)\xi(\tau_1)|g(\tau_1,y_j,w_j)\,dy_j\,dw_j +\int_{\tau_0}^{\tau_1}\int_{\mathbb R^6}|(B\xi)'|g_j\,dy_j\,dw_j\,d\tau\Bigg).
\end{aligned}
\]
Using $C_{B,i}$ and exhausting $J$ by compact subintervals, we obtain
\[
\int_J\int_{\mathbb R^6}\frac{B|\xi|^2}{h+|z|}g_j\,dy_j\,dw_j\,d\tau\lesssim(1+n)C_{B,i}.
\]
The same Cauchy--Schwarz factorization on $|z|\geq h$ gives the fixed-target version of the far-field estimate, with $C_B$ replaced by $C_{B,i}$. On $|z|<h$, Lemma~\ref{lem:rough-support} gives the same near-field bound uniformly in the target trajectory. This proves the fixed-target version of the passage estimate.
\end{remark}

\section{Energy Moment Estimates}\label{sec:energy-moments}

The main goal of this section is to establish weighted bounds for the corrected energy moment over the whole particle cloud and for the associated fourth-order force integral. These cloud-level estimates are the input for the fixed-particle argument in the next section. The proof combines the convex transport estimate in Lemma~\ref{lem:convex-transport} with the passage lemma in Lemma~\ref{lem:passage}.

For the force integrals used below, we write
\[
D_k(\tau)=\int_{\mathbb R^3}m_k(\tau,y)(|K|*m_0)(\tau,y)\,dy.
\]

\begin{proposition}\label{prop:cloud-moment}
Let $\alpha=3/16$. There exist constants $C>0$ and $\gamma>0$, depending only on the initial data, such that for every $T>0$,
\[
Z_T:=1+\sup_{0\leq\tau\leq T}e^{(1+\alpha)\tau}\iint_{\mathbb R^6}H(\tau,y,w)^{5/2}g(\tau,y,w)\,dy\,dw+\int_0^T e^{\gamma\tau}D_4(\tau)\,d\tau\leq C.
\]
\end{proposition}

To prove Proposition~\ref{prop:cloud-moment}, we first establish two estimates for the spatial density and the force integral. We then use the oddness of the Coulomb kernel to exploit the relative motion of particle pairs.

\begin{lemma}\label{lem:density-force}
Let $I\subset[n,n+1]\cap[0,T]$, where $n\geq0$ is an integer. Then
\[
\left\|\int_{\mathbb R^3}g(\tau,y,w)\,dw\right\|_{L^2(I\times\mathbb R^3)}\lesssim e^{3n/4},
\]
and
\[
\int_I\int_{\mathbb R^3}m_1(\tau,y)(|K|*m_0)(\tau,y)\,dy\,d\tau\lesssim e^{7n/16}.
\]
The implicit constants depend only on the initial data.
\end{lemma}

\begin{proof}
Throughout the proof, $I\subset[n,n+1]\cap[0,T]$. By Lemma~\ref{lem:basic-energy},
\[
\sup_I M_2\lesssim e^{-n}.
\]
We first apply the convex transport estimate with weight $1$. To estimate $V_\ell(1)$, we retain the sign of the field term. Let $\Phi(\sigma,y)$ be a smooth convex function with $|\nabla\Phi|\leq1$, and evaluate its derivatives at $(\ell\tau,y_i)$. Along a characteristic,
\[
\frac{d}{d\tau}\left(\partial_\sigma\Phi+\frac{w_i}{\ell}\cdot\nabla_y\Phi\right)=\ell D^2\Phi[(1,w_i/\ell),(1,w_i/\ell)]+\frac{\dot w_i}{\ell}\cdot\nabla_y\Phi.
\]
Integrating this identity in time and in the particle labels by means of \eqref{eq:phase-space-integration}, and substituting $\dot w_i=-w_i+e^{-\tau}E^\sharp(\tau,y_i)$, places a minus sign in front of the field contribution when we solve for the Hessian integral. At each fixed time, the oddness of $K$ gives
\[
\begin{aligned}
\int_{\mathbb R^3}m_0E^\sharp\cdot\nabla_y\Phi\,dy
=\frac12\iint_{\mathbb R^3\times\mathbb R^3}m_0(y_i)m_0(y_j)K(y_i-y_j)\cdot\bigl(\nabla_y\Phi(y_i)-\nabla_y\Phi(y_j)\bigr)\,dy_i\,dy_j\geq0.
\end{aligned}
\]
Indeed, $K(y_i-y_j)$ is a positive multiple of $y_i-y_j$, while the gradient of a convex function is monotone. We may therefore discard this contribution in the upper bound. The endpoint terms are bounded by $2(M_0+\ell^{-1}\sup_I M_1)$, and the damping term by $\ell^{-1}\int_I M_1\,d\tau$. Since $|I|\leq1$ and $M_1\leq M_0^{1/2}M_2^{1/2}$, taking the supremum over $\Phi$ yields
\[
V_\ell(1)\lesssim1+\ell^{-1}e^{-n/2}.
\]
Choose $\ell=e^{-n/2}$. The density estimate in Lemma~\ref{lem:convex-transport} then gives
\begin{equation}\label{eq:density-l2-window}
\left\|\int_{\mathbb R^3}g(\tau,y,w)\,dw\right\|_{L^2(I\times\mathbb R^3)}\lesssim e^n(e^{-n/2})^{1/2}=e^{3n/4}.
\end{equation}

We next estimate $D_1$. The Fourier multiplier of the scalar kernel $|K|$ is a constant multiple of $|\eta|^{-1}$, and the vector multiplier of $K$ has magnitude proportional to $|\eta|^{-1}$; see~\cite{Lieb1983}. Plancherel's identity and \eqref{eq:field-l2} therefore give
\[
\||K|*m_0\|_{L^2_y}\lesssim\|E^\sharp\|_{L^2_y}\lesssim1,\qquad \|\nabla_y(|K|*m_0)\|_{L^2_y}\lesssim\|m_0\|_{L^2_y}.
\]
Sobolev's inequality and interpolation between $L^2_y$ and $L^6_y$ imply
\[
\||K|*m_0\|_{L^4_y}\lesssim\||K|*m_0\|_{L^2_y}^{1/4}\||K|*m_0\|_{L^6_y}^{3/4}\lesssim\|m_0\|_{L^2_y}^{3/4}.
\]
Cauchy--Schwarz in the velocity variable and H\"older's inequality in space now yield
\[
\begin{aligned}
D_1(\tau)&\leq M_2(\tau)^{1/2}\left(\int_{\mathbb R^3}m_0(\tau,y)(|K|*m_0)^2\,dy\right)^{1/2}\\
&\leq M_2(\tau)^{1/2}\|m_0(\tau)\|_{L^2_y}^{1/2}\||K|*m_0\|_{L^4_y}\lesssim e^{-n/2}\|m_0(\tau)\|_{L^2_y}^{5/4}.
\end{aligned}
\]
Applying H\"older's inequality in time and then \eqref{eq:density-l2-window}, we obtain
\begin{equation}
\begin{aligned}
\int_I D_1(\tau)\,d\tau
&\lesssim e^{-n/2}|I|^{3/8}\|m_0\|_{L^2(I\times\mathbb R^3)}^{5/4}\\
&\lesssim e^{-n/2}e^{15n/16}=e^{7n/16}.
\end{aligned}
\label{eq:first-force-window}
\end{equation}
This proves the two estimates and completes the proof.\end{proof}

We next use the oddness of the Coulomb kernel to reveal the relative motion of particle pairs. For two particle labels $i,j$, write $z=y_i-y_j$ and $\xi=w_i-w_j$. Exchanging the two labels in the kinetic contribution of the averaged current gives
\[
\begin{aligned}
 \iint_{\mathbb R^6}\left(\frac{|w_i|^2}{2}\right)^{3/2}J^\sharp(\tau,y_i)g_i\,dy_i\,dw_i  =2^{-5/2}\iint_{\mathbb R^{12}}K(z)\cdot\bigl(|w_i|^3w_j-|w_j|^3w_i\bigr)g_i g_j\,dy_i\,dw_i\,dy_j\,dw_j.
\end{aligned}
\]
Moreover,
\[
\begin{aligned}
\bigl||w_i|^3w_j-|w_j|^3w_i\bigr|
&=|w_i||w_j|\,\bigl||w_i|w_i-|w_j|w_j\bigr|\\
&\leq |w_i||w_j|(|w_i|+|w_j|)|\xi|.
\end{aligned}
\]
Thus the symmetrization produces the relative velocity $\xi$, which is exactly the velocity of the relative position $z$, since $\dot z=\xi$. This is the structure that allows us to apply the passage lemma.

The next lemma controls the variation quantity associated with the particle-pair weight that will be used in the passage estimate.

\begin{lemma}\label{lem:pair-variation}
Let $I\subset[n,n+1]\cap[0,T]$, where $n\geq0$ is an integer, and let $Z_T$ be defined in Proposition~\ref{prop:cloud-moment}. Then
\[
C_{|w_i|^2|w_j|^2}\lesssim Z_T^{1/3}e^{-13n/8}.
\]
The implicit constant depends only on the initial data.
\end{lemma}

\begin{proof}
By the basic energy estimate, Lemma~\ref{lem:density-force}, and the definition of $Z_T$,
\[
\sup_I M_2\lesssim e^{-n},\qquad \sup_I M_5\lesssim Z_Te^{-19n/16},\qquad \int_I D_1(\tau)\,d\tau\lesssim e^{7n/16},\qquad \int_I D_4(\tau)\,d\tau\leq Z_T.
\]
Interpolation between the second and fifth velocity moments gives
\[
\sup_I M_3\leq(\sup_I M_2)^{2/3}(\sup_I M_5)^{1/3}\lesssim Z_T^{1/3}e^{-17n/16}.
\]
Similarly, H\"older's inequality with the common nonnegative weight $(|K|*m_0)g$ gives
\[
\int_I D_2(\tau)\,d\tau\leq\left(\int_I D_1\,d\tau\right)^{2/3}\left(\int_I D_4\,d\tau\right)^{1/3}\lesssim Z_T^{1/3}e^{7n/24}.
\]

For brevity, write $\mathcal B_{ij}=|w_i|^2|w_j|^2$. Since $|\xi|\leq|w_i|+|w_j|$, the endpoint contribution satisfies
\[
\iint_{\mathbb R^{12}}|\mathcal B_{ij}\xi|g_i g_j\,dY\leq2M_2M_3,
\]
where $dY=dy_i\,dw_i\,dy_j\,dw_j$. The damping terms in the two velocity equations contribute $-5\mathcal B_{ij}\xi$, because this expression is homogeneous of degree five in $(w_i,w_j)$.

For the field terms, differentiation gives
\[
\begin{aligned}
\|D_{w_i}(\mathcal B_{ij}\xi)\|&\lesssim |w_i|^2|w_j|^2+|w_i||w_j|^3,\\
\|D_{w_j}(\mathcal B_{ij}\xi)\|&\lesssim |w_i|^2|w_j|^2+|w_i|^3|w_j|.
\end{aligned}
\]
Using $|E^\sharp|\leq |K|*m_0$ and integrating first in the particle variable not carrying the field, the two field contributions are bounded by $M_2D_2+M_3D_1$. Consequently, the definition of $C_{|w_i|^2|w_j|^2}$ gives
\[
\begin{aligned}
C_{|w_i|^2|w_j|^2}
&\lesssim(\sup_I M_2)(\sup_I M_3)+e^{-n}\left[(\sup_I M_2)\int_I D_2(\tau)\,d\tau+(\sup_I M_3)\int_I D_1(\tau)\,d\tau\right]\\
&\lesssim Z_T^{1/3}\left(e^{-33n/16}+e^{-41n/24}+e^{-13n/8}\right)\\
&\lesssim Z_T^{1/3}e^{-13n/8}.
\end{aligned}
\]
This proves the estimate and completes the proof.\end{proof}

\begin{proof}[Proof of Proposition~\ref{prop:cloud-moment}]
Along a characteristic, write $H_i=H(\tau,y_i,w_i)$. Multiplying \eqref{eq:single-energy-2} by $\frac52H_i^{3/2}$ gives
\begin{equation}\label{eq:fifth-energy}
\frac{d}{d\tau}H_i^{5/2}+\frac52H_i^{5/2}+\frac54|w_i|^2H_i^{3/2}=\frac52e^{-\tau}H_i^{3/2}J^\sharp(\tau,y_i).
\end{equation}
Fix $T>0$, and temporarily take $0<\gamma\leq1$ in the definition of $Z_T$. The classical solution has bounded phase-space support on every finite time interval, so $Z_T$ is finite. All constants below are independent of this provisional choice of $\gamma$. Throughout the proof, $I$ denotes an arbitrary subinterval of $[n,n+1]\cap[0,T]$.

\medskip
\noindent\emph{Step 1. The potential contribution.}
We first separate the potential part of $H^{3/2}$. Since $H=|w|^2/2+U^\sharp$ and $U^\sharp\geq0$,
\[
0\leq H^{3/2}-\left(\frac{|w|^2}{2}\right)^{3/2}\lesssim U^\sharp|w|+(U^\sharp)^{3/2}.
\]
Recall that $J^\sharp(\tau,y_i)=\iint K(y_i-y_j)\cdot w_jg_j\,dy_j\,dw_j$. Taking the absolute value of the kernel, using the evenness of $|K|$, and exchanging labels when necessary, we obtain
\[
\iint_{\mathbb R^6}|J^\sharp|g\,dy\,dw\leq D_1,\qquad \iint_{\mathbb R^6}|w||J^\sharp|g\,dy\,dw\leq D_2.
\]
By Lemma~\ref{lem:basic-energy}, $\sup_I\|U^\sharp\|_{L^\infty_y}\lesssim e^{-n/2}$. Lemma~\ref{lem:density-force} and the interpolation estimate for $D_2$ in the proof of Lemma~\ref{lem:pair-variation} therefore give
\begin{equation}\label{eq:potential-work-error}
\begin{aligned}
&\int_I e^{\alpha\tau}\iint_{\mathbb R^6}\left[H^{3/2}-\left(\frac{|w|^2}{2}\right)^{3/2}\right]|J^\sharp|g\,dy\,dw\,d\tau\\
&\lesssim e^{\alpha n}\left(e^{-n/2}\int_I D_2\,d\tau+e^{-3n/4}\int_I D_1\,d\tau\right)\\
&\lesssim Z_T^{1/3}e^{-n/48}+e^{-n/8}.
\end{aligned}
\end{equation}

\medskip
\noindent\emph{Step 2. Symmetrization and the choice of weights.}
For the kinetic contribution, the label-exchange identity established above gives a relative-velocity factor. We apply Lemma~\ref{lem:passage} with the concrete weights
\[
B=|w_i|^2|w_j|^2,\qquad W=e^{\alpha n}|w_i||w_j|(|w_i|+|w_j|)|\xi|+|w_i|^4.
\]
The first term of $W$ controls the symmetrized kinetic contribution, since $e^{\alpha\tau}\leq e^\alpha e^{\alpha n}$ on $I$. The second term, integrated against $|K(z)|g_i g_j$, gives $D_4$. By Lemma~\ref{lem:rough-support}, $|w_i|+|w_j|\lesssim e^n$, and hence $W\lesssim e^{5n}$. Moreover, $B\xi$ is smooth along the classical trajectories. The weight $W$ vanishes when $w_i=0$, while the sets $w_j=0$ and $w_j=w_i$ have zero source-velocity measure. Thus the zero-set condition in Lemma~\ref{lem:passage} holds almost everywhere.

\medskip
\noindent\emph{Step 3. The quotient estimate.}
We estimate the two contributions to the quotient in Lemma~\ref{lem:passage} separately. First, the square of the first term of $W$, divided by $B|\xi|^2$, is
\[
e^{2\alpha n}(|w_i|+|w_j|)^2.
\]
After integration in velocity, this requires a bound on $m_2$. We obtain it from Lemma~\ref{lem:convex-transport} with $\omega=|w|^3$, since $\omega^{2/3}=|w|^2$, and choose $\ell=e^{-n/6}$.

Splitting into $|w|\leq\ell$ and $|w|\geq\ell$ gives
\[
|w|^3+\ell^{-1}|w|^4\leq2\ell^{-2}(|w|^5+\ell^3|w|^2),\qquad
|w|^2+\ell^{-1}|w|^3\leq2\ell^{-2}(|w|^4+\ell^3|w|).
\]
The first inequality controls the endpoint and damping terms in the variation estimate of Lemma~\ref{lem:convex-transport}; the second controls the field terms. Indeed, the damping derivatives of $|w|^3$ and $w|w|^3$ are $-3|w|^3$ and $-4w|w|^3$, while $|\nabla_w|w|^3|=3|w|^2$ and $\|D_w(w|w|^3)\|\leq4|w|^3$. Using $|E^\sharp|\leq|K|*m_0$ and $|I|\leq1$, we obtain
\[
V_\ell(|w|^3)\lesssim\ell^{-2}\sup_I(M_5+\ell^3M_2)+e^{-n}\ell^{-2}\int_I(D_4+\ell^3D_1)\,d\tau.
\]
At the chosen scale, $\ell^3=e^{-n/2}$. The definition of $Z_T$, Lemma~\ref{lem:density-force}, and Lemma~\ref{lem:basic-energy} give
\[
\ell^3\int_I D_1\,d\tau\lesssim e^{-n/16},\qquad \sup_I(M_5+\ell^3M_2)\lesssim e^{-n}Z_T.
\]
Since $e^{-n}\ell^{-2}=e^{-2n/3}$, it follows that
\[
V_\ell(|w|^3)\lesssim e^{-2n/3}Z_T.
\]
The density estimate in Lemma~\ref{lem:convex-transport} therefore yields
\begin{equation}\label{eq:m2-l2-window}
\|m_2\|_{L^2(I\times\mathbb R^3)}\lesssim e^n\ell^{1/2}V_\ell(|w|^3)^{2/3}\lesssim Z_T^{2/3}e^{17n/36}.
\end{equation}

For the kernel in Lemma~\ref{lem:passage}, radial integration gives
\begin{equation}\label{eq:passage-kernel-l1}
\|k_h\|_{L^1}=4\pi\int_0^{C_0e^n}\frac{r^2}{(h+r)^3}\,dr\leq4\pi\log(1+C_0e^n/h)\lesssim1+n.
\end{equation}
Young's inequality in space and Cauchy--Schwarz in time, together with Lemma~\ref{lem:density-force} and \eqref{eq:m2-l2-window}, imply
\[\begin{aligned}
\int_I\int_{\mathbb R^3}m_2(k_h*m_0)\,dy\,d\tau
&\leq\|k_h\|_{L^1}\|m_2\|_{L^2(I\times\mathbb R^3)}\|m_0\|_{L^2(I\times\mathbb R^3)}\\
&\lesssim(1+n)Z_T^{2/3}e^{11n/9}.
\end{aligned}\]
For the kinetic contribution, use $(|w_i|+|w_j|)^2\leq2(|w_i|^2+|w_j|^2)$. Since $k_h$ is even, the two terms give the same spatial integral. Thus its contribution to the squared quotient is bounded by
\begin{equation}\label{eq:quotient-kinetic}
C(1+n)Z_T^{2/3}e^{(13/8-1/36)n},\qquad 2\alpha+\frac{11}{9}=\frac{13}{8}-\frac1{36}.
\end{equation}

We next consider the second term $|w_i|^4$ in $W$. Away from the zero sets already considered, $|w_i|^2\leq2(|w_j|^2+|\xi|^2)$ gives
\[
\frac{|w_i|^8}{B|\xi|^2}\leq2|w_i|^4\left(\frac1{|w_j|^2}+\frac1{|w_j-w_i|^2}\right).
\]
Both singularities are in the source velocity $w_j$. For any $\zeta\in\mathbb R^3$ and $r>0$, splitting into $|w_j-\zeta|<r$ and its complement yields
\[
\int_{\mathbb R^3}\frac{g(\tau,y_j,w_j)}{|w_j-\zeta|^2}\,dw_j\leq4\pi r\|g(\tau)\|_{L^\infty_{y,w}}+r^{-2}m_0(\tau,y_j).
\]
Choosing $r=(m_0(\tau,y_j)/\|g(\tau)\|_{L^\infty})^{1/3}$ when $m_0>0$, and noting that the integral vanishes when $m_0=0$, we obtain, since $\|g(\tau)\|_\infty\lesssim e^{3n}$,
\begin{equation}\label{eq:velocity-singular-integral}
\sup_{\zeta\in\mathbb R^3}\int_{\mathbb R^3}\frac{g(\tau,y_j,w_j)}{|w_j-\zeta|^2}\,dw_j\lesssim e^{2n}m_0(\tau,y_j)^{1/3}.
\end{equation}
Apply this bound with $\zeta=0$ and $\zeta=w_i$. Integration in $w_i$ leaves $m_4(\tau,y_i)$. Hölder's inequality in velocity gives $m_4\leq m_5^{2/3}m_2^{1/3}$. Applying Hölder again in $(\tau,y_i,y_j)$ with the nonnegative weight $k_h(y_i-y_j)$, we obtain
\[
\begin{aligned}
&e^{2n}\int_I\iint_{\mathbb R^3\times\mathbb R^3}k_h(y_i-y_j)m_4(y_i)m_0(y_j)^{1/3}\,dy_i\,dy_j\,d\tau\\
&\leq e^{2n}\left(\|k_h\|_{L^1}\int_I M_5(\tau)\,d\tau\right)^{2/3}\left(\int_I\int_{\mathbb R^3}m_2(k_h*m_0)\,dy\,d\tau\right)^{1/3}\\
&\lesssim(1+n)Z_T^{8/9}e^{(13/8-1/108)n}.
\end{aligned}
\]
Indeed, the powers of $Z_T$ and the exponential factor are
\[
\frac23+\frac13\cdot\frac23=\frac89,\qquad 2-\frac23\cdot\frac{19}{16}+\frac13\cdot\frac{11}{9}=\frac{13}{8}-\frac1{108}.
\]
Since $Z_T\geq1$, combining the kinetic and force contributions gives
\begin{equation}\label{eq:quotient-bound}
\int_I\iint_{\mathbb R^{12}}k_h(z)\frac{W^2}{B|\xi|^2}g_i g_j\,dY\,d\tau\lesssim(1+n)Z_Te^{(13/8-1/108)n}.
\end{equation}

\medskip
\noindent\emph{Step 4. The passage estimate and the time summation.}
By Lemma~\ref{lem:pair-variation},
\[
C_B=C_{|w_i|^2|w_j|^2}\lesssim Z_T^{1/3}e^{-13n/8}.
\]
Combining this with \eqref{eq:quotient-bound} in Lemma~\ref{lem:passage}, we obtain
\[
\begin{aligned}
\int_I\iint_{\mathbb R^{12}}W|K(z)|g_i g_j\,dY\,d\tau\lesssim e^{-3n}+(1+n)Z_T^{2/3}e^{-n/216}.
\end{aligned}
\]
The two terms in $W$ control the symmetrized kinetic current and $\int_I D_4(\tau)\,d\tau$, respectively. Combining this estimate with the potential contribution \eqref{eq:potential-work-error}, and choosing $0<c_0<1/216$, gives
\begin{equation}\label{eq:cloud-window}
\int_I e^{\alpha\tau}\left|\iint_{\mathbb R^6}H^{3/2}J^\sharp g\,dy\,dw\right|d\tau+\int_I D_4(\tau)\,d\tau\lesssim e^{-c_0n}Z_T^{2/3}.
\end{equation}
The factor $1+n$ is absorbed by the slightly weaker exponential decay. The absolute value is taken after integration over the particle cloud, preserving the cancellation from exchanging the particle labels.

It remains to close the estimate for $Z_T$. Integrating \eqref{eq:fifth-energy} over the particle distribution by means of \eqref{eq:phase-space-integration}, and multiplying by $e^{(1+\alpha)\tau}$, we use $\frac52-(1+\alpha)=\frac{21}{16}>0$ to discard the remaining nonnegative energy terms. Thus
\[
\frac{d}{d\tau}\left(e^{(1+\alpha)\tau}\iint_{\mathbb R^6}H^{5/2}g\,dy\,dw\right)\leq\frac52e^{\alpha\tau}\left|\iint_{\mathbb R^6}H^{3/2}J^\sharp g\,dy\,dw\right|.
\]
Integrating up to any time in $[0,T]$ and summing \eqref{eq:cloud-window} over the unit intervals, including the last interval if it is shorter, yields
\[
\sup_{0\leq\tau\leq T}e^{(1+\alpha)\tau}\iint_{\mathbb R^6}H^{5/2}g\,dy\,dw\leq C+CZ_T^{2/3}\sum_{n\geq0}e^{-c_0n}\leq C(1+Z_T^{2/3}).
\]
Choose $0<\gamma<c_0$. On $[n,n+1]$, $e^{\gamma\tau}\leq e^\gamma e^{\gamma n}$, so the force estimate in \eqref{eq:cloud-window} gives
\[
\int_0^T e^{\gamma\tau}D_4(\tau)\,d\tau\lesssim Z_T^{2/3}\sum_{n\geq0}e^{-(c_0-\gamma)n}\lesssim Z_T^{2/3}.
\]
Combining the last two estimates, we find $Z_T\leq C(1+Z_T^{2/3})$. Since $Z_T\geq1$, this implies $Z_T^{1/3}\leq2C$, uniformly in $T$. This completes the proof of Proposition~\ref{prop:cloud-moment}.\end{proof}

\section{Estimates Along Individual Characteristics}\label{sec:proof-main}

Proposition~\ref{prop:cloud-moment} controls weighted energy and force integrals after averaging over the whole particle cloud. We now use these estimates to control the current along one prescribed particle trajectory. The key point is that the fixed-target version of the passage lemma in Remark~\ref{rem:passage-fixed-target} requires only source-particle integration, while the cloud estimates from Section~\ref{sec:energy-moments} provide the necessary moment and force bounds uniformly in the chosen target. Once this current estimate is established, a last-crossing argument will give decay of the relative velocity and hence the desired uniform velocity-support bound.

We fix a particle with initial label in $\operatorname{supp}f_0$. All constants below are independent of this choice. Let $I_n=[n,n+1]$.

\begin{proposition}\label{prop:single-current}
There exist constants $c,C>0$, depending only on the initial data, such that for every $n\in\mathbb N_0$, every target characteristic issuing from an initial label $i\in\operatorname{supp}f_0$, and every interval $J\subset I_n$ on which
\[
|w_i(\tau)|\geq(1+n)^{-3/2},\qquad \tau\in J,
\]
one has
\[
\int_J e^{-\tau}|J^\sharp(\tau,y_i(\tau))|\,d\tau\leq Ce^{-cn}.
\]
\end{proposition}

\begin{proof}
Fix $n$, a target characteristic $i$ issuing from an initial label in $\operatorname{supp}f_0$, and an interval $J$ as in the statement. For brevity, we write $C_B$ for the fixed-target variation quantity $C_{B,i}$ defined in Remark~\ref{rem:passage-fixed-target}; this notation is used only in the present proof. Write $(y_j,w_j)$ for a source particle and $g_j=g(\tau,y_j,w_j)$. Throughout the proof, the target follows its characteristic, and all particle integrals are taken only over the source variables. Set
\[
z=y_i-y_j,\qquad \xi=w_i-w_j.
\]
Then $\dot z=\xi$, and
\[
|w_i(\tau)|\geq(1+n)^{-3/2},\qquad \tau\in J.
\]

The current contains the source velocity $w_j$, whereas the field does not. We estimate both by introducing
\[
T_J=e^{-n}\int_J\iint_{\mathbb R^6}
\bigl(|w_j|+e^{-n/6}\bigr)|K(z)|g_j\,dy_j\,dw_j\,d\tau.
\]
Since $e^{-\tau}\leq e^{-n}$ on $J$, the definitions of the current and the field give
\begin{equation}\label{eq:current-field-target}
\int_J e^{-\tau}|J^\sharp(\tau,y_i)|\,d\tau\leq T_J,
\qquad
e^{-n/6}\int_J e^{-\tau}|E^\sharp(\tau,y_i)|\,d\tau\leq T_J.
\end{equation}
The second inequality will control the field produced when the target velocity is differentiated.

On this finite time interval, the source density and velocity support are bounded. The local integrability of $|K|$, together with the finite total mass, therefore ensures that $T_J<\infty$. We shall obtain a bound independent of the target particle and with exponential decay in $n$.

\medskip
\noindent\emph{Step 1. The passage weight and its variation.}
In the fixed-target version of Lemma~\ref{lem:passage}, take
\[
W=|w_j|+e^{-n/6},\qquad
B=\frac{|w_j|^5+e^{-n/2}|w_j|^2}{|w_i|+|w_j|}.
\]
The denominator gives
\[
|B\xi|\leq |w_j|^5+e^{-n/2}|w_j|^2.
\]
Thus the endpoint terms are controlled by the fifth and second moments of the source particles, uniformly in the target velocity. Differentiation in the source velocity will produce the fourth and first force weights estimated in Section~\ref{sec:energy-moments}.

The rough support bound \eqref{eq:rough-support} gives $W\lesssim e^n\lesssim e^{5n}$. Since the target speed is positive on $J$, the quantity $B|\xi|^2$ vanishes only at $w_j=0$ and $w_j=w_i$. These two velocities have zero source measure. We use the quotient convention of Lemma~\ref{lem:passage} on this set.

We first consider the damping contribution to the derivative of $B\xi$. The common damping scales $w_i$ and $w_j$ by the same factor, so $\xi/(|w_i|+|w_j|)$ is unchanged. The derivative of the numerator under damping is
\[
-|w_j|A'(|w_j|),\qquad
A(q)=q^5+e^{-n/2}q^2.
\]
Since $qA'(q)\leq5A(q)$, the damping contribution is bounded by $5A(|w_j|)$. As $|J|\leq1$, the endpoint and damping terms in $C_B$ are bounded by
\[
C\sup_{\tau\in J}\bigl(M_5(\tau)+e^{-n/2}M_2(\tau)\bigr).
\]

For the field terms, we need the two velocity derivatives. When $w_j\neq0$, the derivative of $\xi/(|w_i|+|w_j|)$ in either velocity has operator norm at most
\[
\frac{1}{|w_i|+|w_j|}
+\frac{|\xi|}{(|w_i|+|w_j|)^2}
\leq\frac{2}{|w_i|+|w_j|}.
\]
Consequently,
\[
\begin{aligned}
\|D_{w_j}(B\xi)\|
&\lesssim |w_j|^4+e^{-n/2}|w_j|,\\
\|D_{w_i}(B\xi)\|
&\leq
\frac{2(|w_j|^5+e^{-n/2}|w_j|^2)}
{|w_i|+|w_j|}\\
&\leq
2(1+n)^{3/2}
\bigl(|w_j|^5+e^{-n/2}|w_j|^2\bigr).
\end{aligned}
\]
At $w_j=0$, both derivatives extend continuously by zero, since $|w_i|\geq(1+n)^{-3/2}$ and the numerator vanishes at least quadratically in $w_j$. In particular, $B\xi$ is absolutely continuous along the particle trajectories on $J$.

The source field contributes
\[
e^{-\tau}D_{w_j}(B\xi)E^\sharp(\tau,y_j).
\]
Using $|E^\sharp|\leq |K|*m_0$, its time and source integrals are bounded by
\[
Ce^{-n}\int_J\bigl(D_4(\tau)+e^{-n/2}D_1(\tau)\bigr)\,d\tau.
\]
For the target field, $E^\sharp(\tau,y_i)$ is independent of the source variables. We therefore integrate the source weight first and use \eqref{eq:current-field-target}:
\[
\begin{aligned}
&\int_J e^{-\tau}|E^\sharp(\tau,y_i)|
\iint_{\mathbb R^6}\|D_{w_i}(B\xi)\|g_j\,dy_j\,dw_j\,d\tau\\
&\qquad\lesssim
(1+n)^{3/2}\sup_J(M_5+e^{-n/2}M_2)
\int_J e^{-\tau}|E^\sharp(\tau,y_i)|\,d\tau\\
&\qquad\leq
(1+n)^{3/2}e^{n/6}
\sup_J(M_5+e^{-n/2}M_2)\,T_J.
\end{aligned}
\]
Combining the endpoint, damping, and field contributions gives
\[
\begin{aligned}
C_B\lesssim{}&
\sup_J(M_5+e^{-n/2}M_2)
+e^{-n}\int_J(D_4+e^{-n/2}D_1)\,d\tau\\
&+(1+n)^{3/2}e^{n/6}
\sup_J(M_5+e^{-n/2}M_2)\,T_J.
\end{aligned}
\]

We now insert the estimates established in Section~\ref{sec:energy-moments}. Proposition~\ref{prop:cloud-moment} and \eqref{eq:basic-moment-energy} give
\[
\sup_J M_5\lesssim e^{-19n/16},
\qquad
\sup_J M_2\lesssim e^{-n},
\qquad
\int_J D_4\,d\tau\lesssim e^{-\gamma n}.
\]
Applying \eqref{eq:first-force-window} to $J$ also gives
\[
\int_J D_1\,d\tau\lesssim e^{7n/16}.
\]
Consequently,
\[
\begin{aligned}
e^n\sup_J(M_5+e^{-n/2}M_2)
&\lesssim e^{-3n/16}+e^{-n/2},\\
e^{7n/6}\sup_J(M_5+e^{-n/2}M_2)
&\lesssim e^{-n/48}+e^{-n/3},\\
e^{-n/2}\int_JD_1\,d\tau&\lesssim e^{-n/16}.
\end{aligned}
\]
In particular, the factor $e^{n/6}$ used to estimate the target field still leaves exponential decay. Adding the contributions to $C_B$ yields
\begin{equation}\label{eq:CB-fixed-target}
e^nC_B\lesssim e^{-\gamma n}+e^{-n/16}
+(1+n)^{3/2}e^{-n/48}T_J.
\end{equation}

\medskip
\noindent\emph{Step 2. The quotient estimate.}
We next estimate the quotient in Lemma~\ref{lem:passage}. Write $q=|w_j|$, $p=|w_i|$, and $r=|\xi|$. Since
\[
A(q)=q^2(q^3+e^{-n/2})
\]
and
\[
(q+e^{-n/6})^3\leq4(q^3+e^{-n/2}),
\]
we have, away from the zero sets,
\[
\frac{W^2}{B|\xi|^2}
\leq
\frac{4(p+q)}{q^2(q+e^{-n/6})r^2}.
\]
To separate the two singularities, use
\[
(1+n)^{-3/2}\leq p\leq q+r.
\]
If $q\leq r$, then $p\leq2r$, $p+q\leq3r$, and $r\geq\frac12(1+n)^{-3/2}$. Thus
\[
\frac{p+q}{r^2}\leq6(1+n)^{3/2}.
\]
If $r\leq q$, the same argument gives
\[
\frac{p+q}{q^2}\leq6(1+n)^{3/2},
\qquad q+e^{-n/6}\geq r+e^{-n/6}.
\]
In both cases,
\[
\frac{W^2}{B|\xi|^2}
\lesssim
(1+n)^{3/2}
\left[
\frac{1}{q^2(q+e^{-n/6})}
+\frac{1}{r^2(r+e^{-n/6})}
\right].
\]
The two terms have the same radial profile, centered at $0$ and $w_i$.

The rough support estimate \eqref{eq:rough-support} gives $|w_j|,|w_i|\leq Ce^n$, with $C$ depending only on the initial data. For either center $\zeta\in\{0,w_i\}$, translation and radial integration therefore give
\[
\begin{aligned}
\int_{|w_j|\leq Ce^n}
\frac{dw_j}{|w_j-\zeta|^2(|w_j-\zeta|+e^{-n/6})}
&\leq4\pi\int_0^{2Ce^n}\frac{ds}{s+e^{-n/6}}\\
&=4\pi\log(1+2Ce^{7n/6})
\lesssim1+n.
\end{aligned}
\]
This bound is uniform in the target. Since $\sup_J\|g\|_\infty\lesssim e^{3n}$, the source velocity integral of the quotient against $g$ is at most $Ce^{3n}(1+n)^{5/2}$ at each source position. By \eqref{eq:passage-kernel-l1}, integration in the source position contributes a factor $\|k_h\|_{L^1}\lesssim1+n$.
Since $|J|\leq1$, we conclude that
\begin{equation}\label{eq:quotient-fixed-target}
\int_J\iint_{\mathbb R^6}
k_h(z)\frac{W^2}{B|\xi|^2}g_j\,dy_j\,dw_j\,d\tau
\lesssim e^{3n}(1+n)^{7/2}.
\end{equation}

\medskip
\noindent\emph{Step 3. Closing the current estimate.}
Apply the fixed-target version of Lemma~\ref{lem:passage}. In terms of $T_J$, it reads
\[
e^nT_J
\lesssim e^{-3n}
+\bigl((1+n)C_B\bigr)^{1/2}
\left(
\int_J\iint_{\mathbb R^6}
k_h(z)\frac{W^2}{B|\xi|^2}g_j\,dy_j\,dw_j\,d\tau
\right)^{1/2}.
\]
Multiplying by $e^{-n}$ and squaring changes the near-field error to $e^{-8n}$. By \eqref{eq:quotient-fixed-target},
\[
T_J^2
\lesssim e^{-8n}+e^{-2n}(1+n)C_Be^{3n}(1+n)^{7/2}
=e^{-8n}+(1+n)^{9/2}e^nC_B.
\]
Inserting \eqref{eq:CB-fixed-target}, we obtain
\[
T_J^2
\lesssim e^{-8n}
+(1+n)^{9/2}(e^{-\gamma n}+e^{-n/16})
+(1+n)^6e^{-n/48}T_J.
\]
Fix $b>0$ such that $2b<\min\{\gamma,1/48\}$. Absorbing the polynomial factors into the exponential bounds gives
\[
T_J^2\leq Ce^{-2bn}(1+T_J).
\]
Young's inequality yields
\[
Ce^{-2bn}T_J
\leq\frac12T_J^2+\frac{C^2}{2}e^{-4bn},
\]
and hence
\[
T_J\lesssim e^{-bn}.
\]
The first inequality in \eqref{eq:current-field-target} now gives
\[
\int_J e^{-\tau}|J^\sharp(\tau,y_i(\tau))|\,d\tau
\lesssim e^{-bn}.
\]
All estimates use only the cloud bounds on $[n,n+1]$, the lower bound on the target speed on $J$, and $|J|\leq1$. The constants are therefore uniform in the target and in the endpoints of $J$. This completes the proof of Proposition~\ref{prop:single-current}.
\end{proof}

The current estimate in Proposition~\ref{prop:single-current} is designed for intervals on which the relative velocity stays above a prescribed threshold. We now show that the corrected energy forces precisely this threshold whenever its suitably weighted value is large. A last-crossing argument then allows us to apply Proposition~\ref{prop:single-current} on each unit time interval and sum the resulting current bounds. This yields decay of the relative velocity, after which the uniform bound on the physical velocity follows from the expanding-coordinate equations.

\Needspace{8\baselineskip}
\begin{proof}[Proof of Theorem~\ref{thm:main}]
Fix a particle characteristic and write $H_i(\tau)=H(\tau,y_i(\tau),w_i(\tau))$. Set
\[
\chi(\tau)=\left(1+\frac{\tau}{4}\right)^3.
\]
By the corrected energy identity \eqref{eq:single-energy-2} and $\chi'(\tau)/\chi(\tau)=3/(4+\tau)\leq3/4$, we have
\begin{equation}\label{eq:weighted-energy-last-crossing}
\frac{d}{d\tau}\bigl(\chi(\tau)H_i(\tau)\bigr)
\leq \chi(\tau)e^{-\tau}|J^\sharp(\tau,y_i(\tau))|.
\end{equation}
The potential estimate \eqref{eq:potential-decay-basic} shows that
\[
\sup_{\tau\geq0}\chi(\tau)U^\sharp(\tau,y_i(\tau))<\infty,
\]
uniformly in the initial particle label. Choose $K_0$ larger than this bound by at least one and also larger than the initial values of $\chi H_i$. Whenever $\chi(\tau)H_i(\tau)\geq K_0$, it follows that
\[
\frac12\chi(\tau)|w_i(\tau)|^2\geq1.
\]
For $\tau\in I_n=[n,n+1]$, we have
\[
\chi(\tau)\leq\left(\frac54\right)^3(1+n)^3<2(1+n)^3,
\]
so the preceding inequality implies
\[
|w_i(\tau)|\geq(1+n)^{-3/2}.
\]

Now fix $\tau_1>0$. If $\chi(\tau_1)H_i(\tau_1)\leq K_0$, there is nothing to prove. Otherwise, since the initial weighted energy is below $K_0$, continuity gives a last time $s<\tau_1$ such that $\chi(s)H_i(s)=K_0$. On $[s,\tau_1]$, the weighted energy remains above $K_0$. Thus the speed threshold required in Proposition~\ref{prop:single-current} holds on the intersection of $[s,\tau_1]$ with each unit interval.

Integrating \eqref{eq:weighted-energy-last-crossing} from $s$ to $\tau_1$ and applying Proposition~\ref{prop:single-current} on these intersections, we obtain
\[
\begin{aligned}
\chi(\tau_1)H_i(\tau_1)
&\leq K_0+C\sum_{n=0}^{\infty}(1+n)^3e^{-cn}\leq C(f_0).
\end{aligned}
\]
The bound is uniform in $\tau_1$ and in the chosen particle. Therefore,
\begin{equation}\label{eq:weighted-energy-final}
\chi(\tau)H_i(\tau)\leq C(f_0),\qquad \tau\geq0.
\end{equation}
 
Since $H_i=|w_i|^2/2+U^\sharp(\tau,y_i)$ and $U^\sharp\geq0$, \eqref{eq:weighted-energy-final} gives
\[
|w_i(\tau)|^2\leq 2C(f_0)\chi(\tau)^{-1},
\]
and hence
\[
|w_i(\tau)|\leq C(f_0)(1+\tau)^{-3/2}.
\]
This estimate is uniform over all initial labels in $\operatorname{supp}f_0$. In particular, the relative velocity is integrable in logarithmic time. Since $\dot y_i=w_i$ and the initial spatial support is compact, we obtain
\[
\begin{aligned}
|y_i(\tau)|
&\leq |y_i(0)|+\int_0^\tau|w_i(\sigma)|\,d\sigma\\
&\leq |y_i(0)|+C(f_0)\int_0^\infty(1+\sigma)^{-3/2}\,d\sigma\leq C(f_0).
\end{aligned}
\]
Finally, the inverse change of variables gives $V_i=y_i+w_i$. Both terms are uniformly bounded, so the physical velocity is uniformly bounded for all time and for every particle in the support. This proves Theorem~\ref{thm:main}.\qedhere
\end{proof}

\begin{samepage}

\nopagebreak[4]
\noindent\textbf{Data/Code Availability.} The manuscript has no associated data or code.
\quad\\

\nopagebreak[4]
\noindent\textbf{Conflict of interest.} The author declares that he has no conflict of interest.
\quad\\

\nopagebreak[4]
\noindent\textbf{AI Declaration.} OpenAI's ChatGPT (GPT-5.6 and GPT-6) assisted with drafting and editing the English exposition, organizing the manuscript, bibliographic checks, and exploratory mathematical discussion. The author has reviewed and edited all AI-assisted material and takes full responsibility for the content of this manuscript.
\end{samepage}

\end{document}